\documentclass[12pt]{amsart}

\usepackage{amsmath,amssymb,amsthm,amscd}

\begin{document}

\newcommand{\End}{{\rm{End}\ts}}
\newcommand{\Hom}{{\rm{Hom}}}
\newcommand{\ch}{{\rm{ch}\ts}}
\newcommand{\non}{\nonumber}
\newcommand{\wt}{\widetilde}
\newcommand{\wh}{\widehat}
\newcommand{\ot}{\otimes}
\newcommand{\la}{\lambda}
\newcommand{\La}{\Lambda}
\newcommand{\al}{\alpha}
\newcommand{\be}{\beta}
\newcommand{\ga}{\gamma}
\newcommand{\ka}{\kappa}
\newcommand{\si}{\sigma}
\newcommand{\vp}{\varphi}
\newcommand{\de}{\delta^{}}
\newcommand{\om}{\omega^{}}
\newcommand{\hra}{\hookrightarrow}
\newcommand{\ve}{\varepsilon}
\newcommand{\ts}{\,}
\newcommand{\vac}{\mathbf{1}}
\newcommand{\di}{\partial}
\newcommand{\qin}{q^{-1}}
\newcommand{\tss}{\hspace{1pt}}
\newcommand{\Sr}{ {\rm S}}
\newcommand{\U}{ {\rm U}}
\newcommand{\Y}{ {\rm Y}}
\newcommand{\AAb}{\mathbb{A}\tss}
\newcommand{\CC}{\mathbb{C}\tss}
\newcommand{\QQ}{\mathbb{Q}\tss}
\newcommand{\SSb}{\mathbb{S}\tss}
\newcommand{\ZZ}{\mathbb{Z}\tss}
\newcommand{\Z}{{\rm Z}}
\newcommand{\Ac}{\mathcal{A}}
\newcommand{\Pc}{\mathcal{P}}
\newcommand{\Qc}{\mathcal{Q}}
\newcommand{\Tc}{\mathcal{T}}
\newcommand{\Sc}{\mathcal{S}}
\newcommand{\Bc}{\mathcal{B}}
\newcommand{\Ec}{\mathcal{E}}
\newcommand{\Hc}{\mathcal{H}}
\newcommand{\Uc}{\mathcal{U}}
\newcommand{\Wc}{\mathcal{W}}
\newcommand{\Ar}{{\rm A}}
\newcommand{\Ir}{{\rm I}}
\newcommand{\Zr}{{\rm Z}}
\newcommand{\gl}{\mathfrak{gl}}
\newcommand{\Pf}{{\rm Pf}}
\newcommand{\oa}{\mathfrak{o}}
\newcommand{\spa}{\mathfrak{sp}}
\newcommand{\g}{\mathfrak{g}}
\newcommand{\z}{\mathfrak{z}}
\newcommand{\Zgot}{\mathfrak{Z}}
\newcommand{\p}{\mathfrak{p}}
\newcommand{\sll}{\mathfrak{sl}}
\newcommand{\agot}{\mathfrak{a}}
\newcommand{\qdet}{ {\rm qdet}\ts}
\newcommand{\cdet}{ {\rm cdet}}
\newcommand{\tr}{ {\rm tr}}
\newcommand{\loc}{{\rm loc}}
\newcommand{\Gr}{ {\rm Gr}\tss}
\newcommand{\sgn}{ {\rm sgn}\ts}
\newcommand{\Sym}{\mathfrak S}
\newcommand{\fand}{\quad\text{and}\quad}
\newcommand{\Fand}{\qquad\text{and}\qquad}
\newcommand{\vt}{{\tss|\hspace{-1.5pt}|\tss}}

\renewcommand{\theequation}{\arabic{section}.\arabic{equation}}

\newtheorem{thm}{Theorem}[section]
\newtheorem{lem}[thm]{Lemma}
\newtheorem{prop}[thm]{Proposition}
\newtheorem{cor}[thm]{Corollary}
\newtheorem{conj}[thm]{Conjecture}

\theoremstyle{definition}
\newtheorem{defin}[thm]{Definition}

\theoremstyle{remark}
\newtheorem{remark}[thm]{Remark}
\newtheorem{example}[thm]{Example}

\newcommand{\bth}{\begin{thm}}
\renewcommand{\eth}{\end{thm}}
\newcommand{\bpr}{\begin{prop}}
\newcommand{\epr}{\end{prop}}
\newcommand{\ble}{\begin{lem}}
\newcommand{\ele}{\end{lem}}
\newcommand{\bco}{\begin{cor}}
\newcommand{\eco}{\end{cor}}
\newcommand{\bde}{\begin{defin}}
\newcommand{\ede}{\end{defin}}
\newcommand{\bex}{\begin{example}}
\newcommand{\eex}{\end{example}}
\newcommand{\bre}{\begin{remark}}
\newcommand{\ere}{\end{remark}}
\newcommand{\bcj}{\begin{conj}}
\newcommand{\ecj}{\end{conj}}

\newcommand{\bal}{\begin{aligned}}
\newcommand{\eal}{\end{aligned}}
\newcommand{\beq}{\begin{equation}}
\newcommand{\eeq}{\end{equation}}
\newcommand{\ben}{\begin{equation*}}
\newcommand{\een}{\end{equation*}}

\newcommand{\bpf}{\begin{proof}}
\newcommand{\epf}{\end{proof}}

\def\beql#1{\begin{equation}\label{#1}}

\title[Sugawara operators]
{On higher order Sugawara operators}

\author{A. V. Chervov}
\address{Institute for Theoretical and Experimental Physics,
Moscow 117259, Russia}
\email{chervov@itep.ru}

\author{A. I. Molev}
\address{School of Mathematics and Statistics,
University of Sydney, NSW 2006, Australia}
\email{alexm\hspace{0.09em}@\hspace{0.1em}maths.usyd.edu.au}

\date{} 


\begin{abstract}
The higher Sugawara operators acting
on the Verma modules over
the affine Kac--Moody algebra at the critical level
are related to the higher Hamiltonians of the Gaudin model
due to work of Feigin, Frenkel and Reshetikhin.
An explicit construction of the higher Hamiltonians
in the case of $\gl_n$ was given recently
by Talalaev.
We propose a new approach to these
results from the viewpoint of the vertex
algebra theory by proving directly the formulas
for the higher order Sugawara operators.
The eigenvalues of the operators
in the Wakimoto modules of critical level are
also calculated.\\[1.0em]
Report ITEP-TH-36/08
\end{abstract}

\maketitle


\section{Introduction}
\label{sec:int}
\setcounter{equation}{0}

Let $\g$ be a simple Lie algebra over $\CC$.
The corresponding affine Kac--Moody algebra $\wh\g$
is defined as the central
extension
\beql{km}
\wh\g=\g\tss[t,t^{-1}]\oplus\CC K,\qquad
\g\tss[t,t^{-1}]=\g\ot\CC[t,t^{-1}],
\eeq
where $\CC[t,t^{-1}]$ is the algebra of Laurent
polynomials in $t$.
Consider the Verma module $M(\la)$ over $\wh\g$ at the critical level,
so that the central element $K$ acts
as multiplication by the negative of the
dual Coxeter number $-h^{\vee}$. The Sugawara operators
form a commuting family of $\wh\g$-endomorphisms of $M(\la)$.
Such families of operators were first constructed
by Goodman and Wallach~\cite{gw:ho} for the $A$ series
and independently by Hayashi~\cite{h:so} for the $A,B,C$ series.
These constructions were used in both papers
for a derivation of the character formula for the
irreducible quotient $L(\la)$ of $M(\la)$,
under certain conditions on $\la$, thus proving
the Kac--Kazhdan conjecture. The existence of
families of Sugawara operators for any
simple Lie algebra $\g$
was established by Feigin and Frenkel~\cite{ff:ak}
by providing a description of the center
of the local completion $\U_{-h^{\vee}}(\wh\g)_{\loc}$
of the universal enveloping algebra
of $\wh\g$ at the critical level. This work is based
on the vertex algebra theory and makes use of the Wakimoto
modules over $\wh\g$; see also~\cite{f:wm}.

On the other hand,
due to work of Feigin, Frenkel and Reshetikhin~\cite{ffr:gm}
the central elements of the local completion
$\U_{-h^{\vee}}(\wh\g)_{\loc}$
are related to the higher Hamiltonians of the Gaudin model
describing quantum spin chain. In the case of $\g=\gl_n$,
Talalaev~\cite{t:qg} produced
a remarkably simple explicit construction
of the higher Hamiltonians.
His construction employs a close connection between
the $XXX$ and Gaudin quantum integrable models
which was extended further by
Mukhin, Tarasov and Varchenko~\cite{mtv:be}. In \cite{mtv:be}
a version of the Bethe ansatz was used to produce
some families of
common eigenvectors of commutative subalgebras of the
Yangian $\Y(\gl_n)$ and the
universal enveloping algebra $\U(\gl_n[t])$
in tensor products
of irreducible $\gl_n$-modules (for an introduction
to the Yangians see \cite{m:yc}).

The relationship between the higher Gaudin Hamiltonians and
the local completion of the universal enveloping algebra
was employed in
\cite{ct:qs} to produce central elements of
$\U_{-h^{\vee}}(\wh\gl_n)_{\loc}$ in an explicit form.
However, the  proof given there
is indirect and relies on
\cite{ff:ak}, \cite{r:uh} and \cite{t:qg}.
In this paper we derive this result directly
with the use of the vertex algebra theory.
We construct two complete sets
of higher Sugawara operators for $\gl_n$ which
thus provides explicit formulas for the
singular vectors of the Verma modules $M(\la)$
at the critical level.
The first set is directly related to \cite{t:qg}
and the construction recovers the corresponding results of
\cite{cf:mm}, \cite{ct:qs},
as well as the main result of \cite{t:qg}.
We also reproduce a description of the classical
$\Wc$-algebra associated to $\gl_n$ and
calculate the eigenvalues of the central elements
in the Wakimoto modules.

\section{Preliminaries on vertex algebras}
\label{sec:pva}
\setcounter{equation}{0}

\subsection{Universal affine
vertex algebra}

A vertex algebra $V$ is a vector space with
the additional data $(Y, T, \vac)$, where
the state-field correspondence $Y$ is a map
$$ Y: V \rightarrow \End V[[z,z^{-1}]], $$
the infinitesimal translation $T$ is an operator $T: V \rightarrow V$,
and $\vac$ is a vacuum vector $\vac \in V$. These data
must satisfy certain axioms; see e.g., \cite{dk:fa},
\cite{fb:va}, \cite{k:va}.
For $a \in V$ we write
\ben
Y(a,z) = \sum_{n\in\ZZ} a_{(n)} z^{-n-1}, \quad
a_{(n)} \in \End V.
\een
In particular, for all $a,b \in V$ we must have
$a_{(n)}\cdot b = 0$ for $n \gg 0$.
The span in $\End V$ of all {\it Fourier coefficients\/} $a_{(n)}$
of all vertex operators $Y(a,z)$ is a Lie subalgebra
of $\End V$ which we denote by  $\Uc$.
The commutator of Fourier coefficients
is given by
\beql{commut}
[a_{(m)},b_{(k)}]=\sum_{n\geqslant 0}\binom{m}{n}
\big(a_{(n)}\cdot b\big)_{(m+k-n)}.
\eeq
Note that \eqref{commut} is equivalent to
the formula for
the commutator $[Y(a,z),Y(b,w)]$ implied by
the operator product expansion formula; see e.g.
\cite{dk:fa}, \cite[Chap.~3]{fb:va}.
The {\it center\/} of the vertex algebra $V$ is its commutative
vertex subalgebra spanned by all vectors $b\in V$
such that $a_{(n)}\cdot b=0$ for all $a\in V$ and $n\geqslant 0$.
The following observation will play an important role below.
If $b$ is an element of the center of $V$, then \eqref{commut}
implies $[Y(a,z), Y(b,w)]=0$ for all $a\in V$.
In other words, all Fourier coefficients $b_{(n)}$
belong to the center of the Lie algebra $\Uc$.

In this paper we are interested in the vertex algebras
associated with affine Kac--Moody algebras $\wh\g$
defined in \eqref{km},
where $\g$ is a finite-dimensional simple
(or reductive) Lie algebra over $\CC$.
Following \cite{dk:fa}, define the {\it universal affine\/}
vertex algebra $V(\g)$ as the quotient of the
universal enveloping algebra $\U(\wh\g)$ by the left
ideal generated by $\g[t]$.
The vacuum vector is $1$ and
the translation operator is determined by
\beql{transl}
T:1\mapsto 0, \quad\big[T, K\big]=0\fand
\big[T, X[n]\tss\big]= -n\tss X[n-1],\qquad X\in\g,
\eeq
where we write $X[n]=X\tss t^{n}$.
The state-field correspondence $Y$
is defined by setting $Y(1,z)=\text{id}$,
\ben
Y(K,z)=K,\qquad Y(J^a[-1],z)=J^a(z)=
\sum_{r\in\ZZ} J^a[r]\tss z^{-r-1},
\een
where $J^1,\dots,J^d$ is a basis of $\g$,
and then extending the map to the whole of $V(\g)$
with the use of {\it normal ordering\/}.
Namely, the normally ordered product of fields
\ben
a(z)=\sum_{r\in\ZZ}a_{(r)}z^{-r-1}\Fand
b(w)=\sum_{r\in\ZZ}b_{(r)}w^{-r-1}
\een
is the formal power series
\beql{normor}
:a(z)\tss b(w){:}= a(z)_+\tss b(w)+b(w)\tss a(z)_-,
\eeq
where
\ben
a(z)_+=\sum_{r<0}a_{(r)}z^{-r-1}\Fand
a(z)_-=\sum_{r\geqslant 0}a_{(r)}z^{-r-1}.
\een
This definition extends to an arbitrary number of fields
with the convention that the normal ordering is read
from right to left. Then
\begin{multline}
Y(J^{a_1}[r_1]\dots J^{a_m}[r_m],z)\\
{}=\frac{1}{(-r_1-1)!\ts\dots (-r_m-1)!}
: \di_z^{-r_1-1}J^{a_1}(z)\dots \di_z^{-r_m-1}J^{a_m}(z):
\non
\end{multline}
with $r_1\leqslant\dots\leqslant r_m<0$, and
$a_i\leqslant a_{i+1}$ if $r_i=r_{i+1}$.

The Lie algebra $\U(\wh\g)_{\loc}$
spanned by the Fourier coefficients of the fields $Y(a,z)$
with $a\in V(\g)$
is known as the {\it local completion\/}
of the universal enveloping algebra
$\U(\wh\g)$ see \cite{ff:ak}, \cite[Sec.~3.5]{fb:va}.
More precisely, the Lie algebra $\U(\wh\g)_{\loc}$
is embedded into the completion $\wt\U(\wh\g)$
of $\U(\wh\g)$ in the topology whose
basis of the open neighborhoods of $0$ is formed by
the left ideals generated by $\g\ot t^N\CC[t]$
with $N\geqslant 0$.

\subsection{Segal--Sugawara vectors}

As a vector space, $V(\g)$ can be identified with
the space of polynomials in $K$ with coefficients
in the universal enveloping algebra
$\U(\g_-)$, where $\g_-=\g\ot t^{-1}\CC[t^{-1}]$.
An element $S\in \U(\g_-)\subseteq V(\g)$ is called
a {\it Segal--Sugawara vector\/} if
\beql{ssvec}
\g\tss[t]\ts S\in (K+h^{\vee})\ts V(\g),
\eeq
where $h^{\vee}$ is the dual Coxeter number for $\g$.
A well-known example of such a vector is provided by
the quadratic element
\beql{sclass}
S=\frac12\ts\sum_{a=1}^d J^a[-1]\tss J_a[-1],
\eeq
where $J_1,\dots,J_d$ is the basis of $\g$
dual to $J^1,\dots,J^d$
with respect to the
suitably normalized invariant bilinear form; see e.g.
\cite[Chap.~2]{fb:va}.

The Segal--Sugawara vectors form a
vector subspace of $\U(\g_-)$ which we will denote by
$\z(\wh\g)$. Moreover, $\z(\wh\g)$ is clearly
closed under the multiplication
in $\U(\g_-)$. The structure of the algebra $\z(\wh\g)$
was described in \cite{ff:ak}; see also \cite{f:wm}.
In order to formulate the result,
for any element $S\in \U(\g_-)$ denote by
$\overline S$ its image in the associated graded algebra
$\text{gr}\ts\U(\g_-)\cong \Sr(\g_-)$.
We call the set of elements
\ben
S_1,\dots,S_n\in \U(\g_-),\qquad n=\text{rk}\ts\g
\een
a {\it complete set of Segal--Sugawara vectors\/} if
each $S_l$ satisfies \eqref{ssvec}
and the corresponding elements
$\overline S_1,\dots,\overline S_n$ coincide with
the images of certain algebraically independent
generators of the algebra
of invariants $\Sr(\g)^{\g}$ under the embedding
$\Sr(\g)\hookrightarrow \Sr(\g_-)$ defined by
the assignment $J^a\mapsto J^a[-1]$.
In accordance to \cite{ff:ak},
$\z(\wh\g)$ is the algebra of polynomials in infinitely
many variables,
\beql{genz}
\z(\wh\g)=\CC[T^{\tss r}S_l\ |\ l=1,\dots,n,\ \ r\geqslant 0].
\eeq
A different proof of this result was given in a more recent work
\cite[Theorem~9.6]{f:wm}.

Equivalently, the algebra $\z(\wh\g)$ can be viewed as the
center of the vertex algebra $V_{-h^{\vee}}(\g)$, where
for each $\ka\in\CC$ the {\it universal affine\/}
vertex algebra $V_{\ka}(\g)$ {\it of level\/}
$\ka$ is the quotient
of $V(\g)$ by the ideal $(K-\ka)\tss V(\g)$. The space
$V_{\ka}(\g)$ is also called the {\it vacuum representation
of\/} $\wh\g$ {\it of level\/} $\ka$.
The center of $V_{\ka}(\g)$
is trivial (i.e., coincides with the multiples
of the vacuum vector) for all values
of $\ka$ except for the critical value $\ka=-h^{\vee}$.
If $\ka\ne-h^{\vee}$ then
the vertex algebra $V_{\ka}(\g)$ has the standard
conformal structure determined by the
Segal--Sugawara vector \eqref{sclass}.
Due to the
celebrated {\it Sugawara construction\/},
the Fourier coefficients of the field
\ben
\frac{1}{2\tss(\ka+h^{\vee})}\ts
\sum_{a=1}^d : J^a(z)\tss J_a(z):
\een
generate an action of the Virasoro algebra on
the space $V_{\ka}(\g)$. On the other hand, if $\ka=-h^{\vee}$,
then the Fourier coefficients of the field
\ben
\sum_{a=1}^d : J^a(z)\tss J_a(z):
\een
are the {\it Sugawara operators\/}, they
commute with the action of
$\wh\g$ on the Verma module $M(\la)$ at the critical level.
The existence of higher order Sugawara operators
on the Verma modules
was established in \cite{gw:ho} and \cite{h:so}, and
certain complete (or fundamental) sets
of such operators were produced
for types $A,B,C$; see also \cite{bbss:ev}.

Given a complete set of Segal--Sugawara vectors
$S_l\in V(\g)$ with $l=1,\dots,n$,
the corresponding fields
\beql{fiese}
S_l(z)=Y(S_l,z)=\sum_{r\in\ZZ} S_{l,(r)}\ts z^{-r-1}
\eeq
will form a {\it complete
set of Sugawara fields\/}.
This terminology
can be shown to be consistent
with \cite{gw:ho} and \cite{h:so}.
The ({\it higher\/}) {\it Sugawara operators\/} are the Fourier
coefficients $S_{l,(r)}\in\U(\wh\g)_{\loc}$
of the fields \eqref{fiese}. If the highest weight $\la$
of the Verma module $M(\la)$ over $\wh\g$ is of
critical level, then $K+h^{\vee}$ is the zero operator
on $M(\la)$. Due to \eqref{commut} and
\eqref{ssvec}, the Sugawara operators
form a commuting family
of $\wh\g$-endomorphisms of $M(\la)$.
Moreover, an explicit
construction of a complete set of Sugawara fields \eqref{fiese}
yields an explicit description of all singular vectors
of the generic Verma modules $M(\la)$ as polynomials
in the Sugawara operators $S_{l,(r)}$ with $r<0$
applied to the highest vector
of $M(\la)$; see \cite[Proposition~4.10]{h:so}
for a precise statement.

The local completion {\it at the critical level\/} is the quotient
$\U_{-h^{\vee}}(\wh\g)_{\loc}$ of $\U(\wh\g)_{\loc}$ by the ideal
generated by $K+h^{\vee}$.
The Sugawara operators can be viewed as the elements of the
center of $\U_{-h^{\vee}}(\wh\g)_{\loc}$. In accordance
to the description
of the center which was originally given
in \cite{ff:ak}, all central elements
are obtained from $\z(\wh\g)$
by the state-field correspondence.

\section{Segal--Sugawara vectors for $\gl_n$}
\label{sec:vr}
\setcounter{equation}{0}

We let $e_{ij}$ denote the standard basis elements of
the Lie algebra $\gl_n$ which we equip with the
invariant symmetric bilinear form
\ben
(X,\ts Y)=\tr\tss(X\tss Y)-\frac{1}{n}\ts\tr\tss X
\ts\tr\tss Y,\qquad X, Y\in\gl_n.
\een
The element $e_{11}+\dots+e_{nn}$ spans the kernel
of the form, therefore it defines a non-degenerate
invariant symmetric bilinear form
on the Lie algebra $\sll_n$.
The elements $e_{ij}t^r$ of $\gl_n[t,t^{-1}]$ with $r\in\ZZ$ will
be denoted by $e_{ij}[r]$.
Then the affine Kac--Moody algebra $\wh\gl_n=\gl_n[t,t^{-1}]\oplus\CC K$
has the commutation relations
\beql{commrel}
\big[e_{ij}[r],e_{kl}[s\tss]\tss\big]
=\de_{kj}\ts e_{i\tss l}[r+s\tss]
-\de_{i\tss l}\ts e_{kj}[r+s\tss]
+K\Big(\de_{kj}\tss\de_{i\tss l}
-\frac{\de_{ij}\tss\de_{kl}}{n}\Big)\ts r\tss\de_{r,-s},
\eeq
and the element $K$ is central. The normalization of the form
is chosen in such a way that the
critical level $-n$ coincides with
the negative of the dual Coxeter number for $\sll_n$.
We will also need
the extended Lie algebra $\wh\gl_n\oplus \CC \tau$,
where for the element $\tau$ we
have the relations
\beql{tauco}
\big[\tau,e_{ij}[r]\tss\big]=-r\ts e_{ij}[r-1],\qquad
\big[\tau,K\big]=0.
\eeq

We identify the universal affine vertex algebra
$V(\gl_n)$ with the vector space of polynomials in $K$
with coefficients in
$\U(\g_-)$  with
$\g_-=t^{-1}\gl_n[t^{-1}]$.
We will also consider the quotient
of $\U(\wh\gl_n\oplus \CC \tau)$ by the left ideal
generated by $\gl_n[t]$.
As a vector space, this quotient will be
identified with
$
V(\gl_n)\ot \CC[\tau],
$
which may be viewed as
the space of polynomials in $\tau$
with coefficients in $V(\gl_n)$.

For an arbitrary $n\times n$
matrix $A=[a_{ij}]$ with entries in a ring
we define its {\it column-determinant\/} $\cdet\ts A$ by
the formula
\beql{cdetdef}
\cdet\ts A=\sum_{\si}\sgn\si\cdot a_{\si(1)1}\dots
a_{\si(n)n},
\eeq
summed over all permutations $\si$ of the set $\{1,\dots,n\}$.

Consider the $n\times n$ matrix
$\tau+E[-1]$
with entries in $\U(\g_-\oplus \CC \tau)$ given by
\ben
\tau+E[-1]=\left[\begin{matrix}
\tau+e_{11}[-1]&e_{12}[-1]&\dots&e_{1n}[-1]\\
e_{21}[-1]&\tau+e_{22}[-1]&\dots&e_{2n}[-1]\\
\vdots&\vdots& \ddots&\vdots     \\
e_{n1}[-1]&e_{n2}[-1]&\dots&\tau+e_{nn}[-1]
                \end{matrix}\right].
\een
Its column-determinant $\cdet\big(\tau+E[-1]\big)$ is an element
of $\U(\g_-\oplus \CC \tau)$ which we also regard as an element of
$V(\gl_n)\ot \CC[\tau]$.

Our main result is a direct proof of the following theorem.

\bth\label{thm:sugaw}
The coefficients $S_1,\dots,S_n$ of the polynomial
\beql{defsl}
\cdet\big(\tau+E[-1]\big)=\tau^n+S_1\tss\tau^{n-1}+\dots+
S_{n-1}\tss\tau+S_n
\eeq
form a complete set of
Segal--Sugawara vectors in $V(\gl_n)$. Hence,
$\z(\wh\gl_n)$
is the algebra of polynomials,
\ben
\z(\wh\gl_n)=\CC[T^{\tss r}S_l\ |\
l=1,\dots,n;\ r\geqslant 0\tss].
\een
\eth

We will prove the theorem in Section~\ref{sec:pro}.
Here we consider
some examples and corollaries, and derive
the second construction
of a complete set of Segal--Sugawara vectors; see
Theorem~\ref{thm:traces} below.

First we point out that regarding the Lie
algebra $\sll_n$ as the quotient of $\gl_n$
by the relation $e_{11}+\dots+e_{nn}=0$,
we obtain the respective complete set of
Segal--Sugawara vectors in $V(\sll_n)$.
In particular, the vector $S_1$ vanishes, while
$-S_2$ coincides with the vector \eqref{sclass}.

In what follows, the usual vertical line
notation for determinants
will always be used for column-determinants.

\bex\label{ex:neqtwothr}
If $n=2$ then
\ben
\bal
S_1&=e_{11}[-1]+e_{22}[-1],\\
S_2&=
\left|\ts
\begin{matrix}
e_{11}[-1]&e_{12}[-1]\\
e_{21}[-1]&e_{22}[-1]
\end{matrix}
\ts\right|+e_{22}[-2]
=e_{11}[-1]\ts e_{22}[-1]-e_{21}[-1]\ts e_{12}[-1]+e_{22}[-2].
\eal
\een
If $n=3$ then
\ben
\bal
S_1&=e_{11}[-1]+e_{22}[-1]+e_{33}[-1],\\
S_2&=\left|\ts
\begin{matrix}
e_{11}[-1]&e_{12}[-1]\\
e_{21}[-1]&e_{22}[-1]
\end{matrix}
\ts\right|
+\left|\ts
\begin{matrix}
e_{11}[-1]&e_{13}[-1]\\
e_{31}[-1]&e_{33}[-1]
\end{matrix}
\ts\right|
+\left|\ts
\begin{matrix}
e_{22}[-1]&e_{23}[-1]\\
e_{32}[-1]&e_{33}[-1]
\end{matrix}
\ts\right|\\
{}&{}+\ts e_{22}[-2]+2\ts e_{33}[-2]
\eal
\een
and
\ben
\bal
S_3&=\left|\ts
\begin{matrix}
e_{11}[-1]&e_{12}[-1]&e_{13}[-1]\\
e_{21}[-1]&e_{22}[-1]&e_{23}[-1]\\
e_{31}[-1]&e_{32}[-1]&e_{33}[-1]
\end{matrix}
\ts\right|+2\ts e_{33}[-3]\\
{}&{}+\left|\ts
\begin{matrix}
e_{11}[-1]&e_{13}[-2]\\
e_{31}[-1]&e_{33}[-2]
\end{matrix}
\ts\right|
+\left|\ts
\begin{matrix}
e_{22}[-1]&e_{23}[-2]\\
e_{32}[-1]&e_{33}[-2]
\end{matrix}
\ts\right|
+\left|\ts
\begin{matrix}
e_{22}[-2]&e_{23}[-1]\\
e_{32}[-2]&e_{33}[-1]
\end{matrix}
\ts\right|.
\eal
\een
\vskip-0.8\baselineskip
\qed
\eex

Now we use Theorem~\ref{thm:sugaw} to obtain
the second complete set of
Segal--Sugawara vectors in $V(\gl_n)$.

\bth\label{thm:traces}
For any $k\geqslant 0$
all coefficients $T_{kl}$ in the expansion
\ben
\tr(\tau+E[-1])^k=
T_{k\tss 0}\ts\tau^{k}+T_{k1}\ts\tau^{k-1}
+\dots+T_{kk}
\een
are Segal--Sugawara vectors in $V(\gl_n)$.
Moreover, the elements $T_{11},\dots,T_{nn}$ form
a complete set of Segal--Sugawara vectors. Hence,
$\z(\wh\gl_n)$
is the algebra of polynomials,
\ben
\z(\wh\gl_n)=\CC[T^{\tss r}T_{l\tss l}\ |\
l=1,\dots,n;\ r\geqslant 0\tss].
\een
\eth

\bpf
Observe that if we replace
$\tau$ by $u+\tau$,
where $u$ is a complex variable,
then relations \eqref{tauco}
will still hold.
Therefore, the coefficients in the expansion
of the polynomial $\cdet(u+\tau+E[-1])$
in the powers of $u$ and $\tau$ are
Segal--Sugawara vectors in $V(\gl_n)$.
Hence, the first part of the theorem
is implied by the
identity
\beql{newton}
\di_u\ts\cdet(u+\tau+E[-1])=\cdet(u+\tau+E[-1])\ts
\sum_{k=0}^{\infty} (-1)^k\ts u^{-k-1}\ts \tr(\tau+E[-1])^k,
\eeq
since it provides an expression for the vectors $T_{kl}$
in terms of the $S_l$ regarded as elements of $\U(\g_-)$.
The identity \eqref{newton} can be viewed as a noncommutative analogue
of the Liouville formula (as well as
the Newton identities).
It follows from \cite[Theorem~4]{cf:mm},
due to Lemma~\ref{lem:tauema}
below.\footnote{We thank Pavel Pyatov for pointing out
that Theorem~4 in \cite{cf:mm} can also be proved
with the use of $R$-matrix arguments; cf. \cite{iop:qm},
\cite[Chap.~7]{m:yc}.}

In order to prove the second part, note that
the elements $\overline T_{11},\dots,\overline T_{nn}$
coincide with the images of the traces
of powers $\tr\ts E^k$ of the matrix $E=[e_{ij}]$
with $k=1,\dots,n$
under the embedding
$\Sr(\gl_n)\hookrightarrow \Sr(\g_-)$ defined by
$e_{ij}\mapsto e_{ij}[-1]$.
However, the elements $\tr\ts E,\dots,\tr\ts E^n$
are algebraically independent and
generate the algebra
of invariants $\Sr(\gl_n)^{\gl_n}$.
\epf

\bex\label{ex:traces}
We have
\ben
\bal
T_{10}&=n,\qquad T_{11}=\tr\tss E[-1]\\
T_{20}&=n,\qquad T_{21}=2\ts\tr\tss E[-1],\qquad
T_{22}=\tr\tss E[-1]^2+\tr\tss E[-2],\\
T_{30}&=n,\qquad T_{31}=3\ts\tr\tss E[-1],\qquad
T_{32}=3\ts\tr\tss E[-1]^2+3\ts\tr\tss E[-2]
\eal
\een
and
\ben
T_{33}=\tr\tss E[-1]^3+2\ts\tr\tss E[-1]\tss E[-2]+
\tr\tss E[-2]\tss E[-1]+2\ts\tr\tss E[-3].
\een
\vskip-1.2\baselineskip
\qed
\eex

\bre\label{rem:oth}
Some other families of Segal--Sugawara vectors
can be constructed by using properties of Manin matrices.
In particular, the quantum MacMahon
master theorem \cite{glz:qm}
leads to a construction of permanent-type vectors; see
also \cite{cf:mm}.
\qed
\ere

\section{Proof of Theorem~\ref{thm:sugaw}}
\label{sec:pro}
\setcounter{equation}{0}

We start by establishing some general properties
of column-determinants.
First we note that the column-determinant $\cdet\ts A$
of a matrix $A$ over an arbitrary ring
changes sign if two rows of the matrix $A$ are swapped.
In particular, $\cdet\ts A=0$ if $A$ has two identical rows.

The next
lemma is immediate from the definition
of the column-determinant.

\ble\label{lem:comm}
Let $A=[a_{ij}]$ be an arbitrary $n\times n$ matrix with entries
in a ring and $b$ an element of the ring.
Then
\ben
\big[\tss b,\cdet\ts A\big]=\sum_{i=1}^n\ts
\left|\ts\begin{matrix}
a_{11}&\cdots&[b,a_{1i}]&\cdots&a_{1n}\\
a_{21}&\cdots&[b,a_{2i}]&\cdots&a_{2n}\\
\cdots&\cdots&\cdots&\cdots&\cdots\\
a_{n1}&\cdots&[b,a_{ni}]&\cdots&a_{nn}
\end{matrix}\ts\right|.
\een
\ele

Lemma~\ref{lem:comm} implies one more property
of column-determinants.

\ble\label{lem:comcol}
Let $A=[a_{ij}]$ be an arbitrary $n\times n$ matrix with entries
in a ring and $b$ an element of the ring. Replace column $j$
of $A$ by the column whose all entries are zero except
for the $i$-th entry equal to $b$.
Then the column-determinant of this matrix can be
written as
\ben
\bal
\left|\ts\begin{matrix}
a_{11}&\cdots&0&\cdots&a_{1n}\\
\cdots&\cdots&\cdots&\cdots&\cdots\\
a_{i1}&\cdots&b&\cdots&a_{in}\\
\cdots&\cdots&\cdots&\cdots&\cdots\\
a_{n1}&\cdots&0&\cdots&a_{nn}
\end{matrix}\ts\right|
&{}=(-1)^{n-j}\left|\ts\begin{matrix}
a_{11}&\cdots&\cdots&a_{1n}&0\\
\cdots&\cdots&\cdots&\cdots&\cdots\\
a_{i1}&\cdots&\cdots&a_{in}&b\\
\cdots&\cdots&\cdots&\cdots&\cdots\\
a_{n1}&\cdots&\cdots&a_{nn}&0
\end{matrix}\ts\right|\\
&{}+(-1)^{i+j}\sum_{k=j+1}^n
\left|\ts\begin{matrix}
a_{11}&\cdots&[b,a_{1k}]&\cdots&a_{1n}\\
a_{21}&\cdots&[b,a_{2k}]&\cdots&a_{2n}\\
\cdots&\cdots&\cdots&\cdots&\cdots\\
a_{n1}&\cdots&[b,a_{nk}]&\cdots&a_{nn}
\end{matrix}\ts\right|,
\eal
\een
where the first determinant on the right hand side
is obtained by moving column $j$ to the last position, while
row $i$ and column $j$ in the determinants occurring in the sum
are deleted and the commutators $[b,a_{i\tss k}]$ occur in column $k-1$.
\qed
\ele

Now we recall some properties of a class of matrices
introduced by Manin~\cite{m:qg}.
Following \cite{cf:mm} we call a matrix $A=[a_{ij}]$
with entries in a ring a
{\it Manin matrix\/} if
\ben
a_{ij}\ts a_{kl}-a_{kl}\ts a_{ij}=a_{kj}\ts a_{il}-a_{il}\ts a_{kj}
\qquad\text{for all possible}\quad i,j,k,l.
\een
Such matrices are also known in the literature as
{\it right-quantum matrices\/} (with $q=1$); see \cite{glz:qm}.
It is straightforward to verify that
the column-determinant of a (square) Manin matrix
will change sign if two columns are swapped.
In particular, if a Manin matrix has two identical
columns, then its column-determinant is zero.
The next observation will be useful in the calculations below.

\ble\label{lem:tauema}
The matrix $\tau+E[-1]$
with entries in $\U(\g_-\oplus \CC \tau)$
is a Manin matrix.
\qed
\ele

We begin proving Theorem~\ref{thm:sugaw}
by verifying that the elements $S_1,\dots,S_n$
satisfy \eqref{ssvec}. Since
\ben
\big[e_{ij}[0\tss],\tau\big]=0\Fand
\big[e_{nn}[1],\tau^k\big]=k\ts \tau^{k-1} e_{nn}[0\tss],
\een
it will be sufficient to check that for all $i,j$
\beql{annih}
e_{ij}[0]\ts \cdet\big(\tau+E[-1]\big)=0
\fand
e_{nn}[1]\ts \cdet\big(\tau+E[-1]\big)\in
(K+n)\tss V(\gl_n)
\eeq
in the $\wh\gl_n$-module $V(\gl_n)\ot \CC[\tau]$.
The first
relation in \eqref{annih} is analogous to
the well-known property of
the Capelli determinant $\cdet[\tss\de_{ij}(u-i+1)+e_{ij}]$;
this is a polynomial in $u$ with coefficients in the
center of the universal enveloping algebra $\U(\gl_n)$.
This property can be verified by a direct argument, and we
argue in a similar way to prove the relation.
By Lemma~\ref{lem:comm},
the polynomial $e_{ij}[0]\ts \cdet\big(\tau+E[-1]\big)$
equals
\beql{sumdet}
\left|\ts\ts\begin{matrix}
\cdots&-e_{1j}[-1]&\cdots\\
\cdots&\cdots&\cdots\\
\cdots&e_{ii}[-1]-e_{jj}[-1]&\cdots\\
\cdots&\cdots&\cdots\\
\cdots&-e_{nj}[-1]&\cdots
\end{matrix}\ts\ts\right|
+\sum_{k\ne i}\ts
\left|\ts\ts\begin{matrix}
\cdots&0&\cdots\\
\cdots&\cdots&\cdots\\
\cdots&e_{ik}[-1]&\cdots\\
\cdots&\cdots&\cdots\\
\cdots&0&\cdots
\end{matrix}\ts\ts\right|,
\eeq
where the dots indicate the same
entries as in the column-determinant
of the matrix
$\tau+E[-1]$,
except for column $i$ in the first
determinant and column $k$ in the $k$-th term in the sum;
the entries shown in the middle belong to row $j$
of each determinant. The first determinant can be written
as the difference of the determinants
\beql{diffdet}
\left|\ts\ts\begin{matrix}
\cdots&0&\cdots\\
\cdots&\cdots&\cdots\\
\cdots&\tau+e_{ii}[-1]&\cdots\\
\cdots&\cdots&\cdots\\
\cdots&0&\cdots
\end{matrix}\ts\ts\right|
-\left|\ts\ts\begin{matrix}
\cdots&e_{1j}[-1]&\cdots\\
\cdots&\cdots&\cdots\\
\cdots&\tau+e_{jj}[-1]&\cdots\\
\cdots&\cdots&\cdots\\
\cdots&e_{nj}[-1]&\cdots
                \end{matrix}\ts\ts\right|.
\eeq
Now, if $i=j$ then the second determinant here equals
$\cdet(\tau+E[-1])$, while the first
determinant together with the sum over $k$ in \eqref{sumdet}
equals $\cdet(\tau+E[-1])$ by the analogue
of the $j$-th row expansion formula for column-determinants.
Hence the first relation
in \eqref{annih} holds in this case.

If $i\ne j$, then the second
determinant in \eqref{diffdet} is obtained from
$\cdet(\tau+E[-1])$ by replacing column $i$
with column $j$.
Similarly, by the row expansion formula
for column-determinants,
the first determinant in \eqref{diffdet}
together with the sum over $k$ in \eqref{sumdet}
is obtained from
$\cdet(\tau+E[-1])$ by replacing row $j$
with row $i$. Both determinants are zero
due to Lemma~\ref{lem:tauema} and the properties
of column-determinants.

Now we turn to the second relation in \eqref{annih}.
By \eqref{commrel}
and \eqref{tauco} we have
\ben
\big[e_{nn}[1],\tau+e_{ii}[-1]\ts\big]=e_{nn}[0]-K'
\Fand
\big[e_{nn}[1],e_{ni}[-1]\ts\big]=e_{ni}[0],
\een
where $i\ne n$ and we have put $K'=K/n$, while
\ben
\big[e_{nn}[1],e_{in}[-1]\ts\big]=-e_{in}[0]
\Fand
\big[e_{nn}[1],\tau+e_{nn}[-1]\ts\big]=(n-1)\ts K'+e_{nn}[0].
\een
Therefore, by Lemma~\ref{lem:comm}, the polynomial
$e_{nn}[1]\ts \cdet\big(\tau+E[-1]\big)$
equals
\beql{sumdetone}
\sum_{i=1}^{n-1}\ts
\left|\ts\ts\begin{matrix}
\cdots&0&\cdots\\
\cdots&\cdots&\cdots\\
\cdots&e_{nn}[0\tss]-K'&\cdots\\
\cdots&0&\cdots\\
\cdots&\cdots&\cdots\\
\cdots&e_{ni}[0\tss]&\cdots
\end{matrix}\ts\ts\right|+
\left|\ts\ts\begin{matrix}
\cdots&-e_{1n}[0\tss]\\
\cdots&\cdots\\
\cdots&-e_{in}[0\tss]\\
\cdots&\cdots\\
\cdots&-e_{n-1,n}[0\tss]\\
\cdots&(n-1)\tss K'+e_{nn}[0\tss]
\end{matrix}\ts\ts\right|,
\eeq
where the dots replace the entries
of the determinant $\cdet(\tau+E[-1])$, except for
column $i$ in the $i$-th summand and the last column
in the last summand. Since each element $e_{in}[0\tss]$
annihilates the vacuum vector of $V(\gl_n)$,
the last summand can be written as
$(n-1)\tss K'\ts\big|\tau+E[-1]\big|_{nn}$, where by
$\big|\tau+E[-1]\big|_{ii}$
we denote the column-determinant of the matrix obtained from
$\tau+E[-1]$ by deleting row and column $i$.
Using this notation, we can represent \eqref{sumdetone}
in the form
\beql{sumdetsi}
\sum_{i=1}^{n-1}\ts
\left|\ts\ts\begin{matrix}
\cdots&0&\cdots\\
\cdots&\cdots&\cdots\\
\cdots&e_{nn}[0\tss]&\cdots\\
\cdots&\cdots&\cdots\\
\cdots&e_{ni}[0\tss]&\cdots
\end{matrix}\ts\ts\right|-K'\sum_{i=1}^{n-1}
\big|\tau+E[-1]\big|_{ii}+
(n-1)\tss K'\ts\big|\tau+E[-1]\big|_{nn}.
\eeq
Our next step is to prove that the following key relation
holds for each value of $i=1,\dots,n-1$:
\beql{ifix}
\left|\ts\ts\begin{matrix}
\cdots&0&\cdots\\
\cdots&\cdots&\cdots\\
\cdots&e_{nn}[0\tss]&\cdots\\
\cdots&\cdots&\cdots\\
\cdots&e_{ni}[0\tss]&\cdots
\end{matrix}\ts\ts\right|=
\big|\tau+E[-1]\big|_{nn}-\big|\tau+E[-1]\big|_{ii}.
\eeq
The determinant on the left hand side equals the sum
of two determinants obtained by replacing
$e_{ni}[0\tss]$ or $e_{nn}[0\tss]$ by $0$, respectively.
Now we apply Lemma~\ref{lem:comcol} to each of these
determinants. For the first one we have
\begin{align}
\label{firde}
&\left|\ts\ts\begin{matrix}
\cdots&0&\cdots\\
\cdots&\cdots&\cdots\\
\cdots&e_{nn}[0\tss]&\cdots\\
\cdots&\cdots&\cdots\\
\cdots&0&\cdots
\end{matrix}\ts\ts\right|=(-1)^{n-i}
\left|\ts\ts\begin{matrix}
\cdots&\cdots&0\\
\cdots&\cdots&\cdots\\
\cdots&\cdots&e_{nn}[0\tss]\\
\cdots&\cdots&\cdots\\
\cdots&\cdots&0
\end{matrix}\ts\ts\right|\\
&{}+\sum_{k=i+1}^{n-1}\left|\ts\ts\begin{matrix}
\cdots&0&\cdots\\
\cdots&0&\cdots\\
\cdots&\cdots&\cdots\\
\cdots&e_{nk}[-1]&\cdots
\end{matrix}\ts\ts\right|
+\left|\ts\ts\begin{matrix}
\cdots&\cdots&-e_{1n}[-1]\\
\cdots&\cdots&\cdots\\
\cdots&\cdots&-e_{n-1,n}[-1]\\
\cdots&\cdots&0
\end{matrix}\ts\ts\right|,
\non
\end{align}
where the determinants in the second line are of
the size $(n-1)\times(n-1)$ with row and column $i$ deleted.
Now add and subtract the
following $(n-1)\times(n-1)$ determinant
with row and column $i$ deleted
\ben
\left|\ts\ts\begin{matrix}
\cdots&\cdots&0\\
\cdots&\cdots&\cdots\\
\cdots&\cdots&\cdots\\
\cdots&\cdots&-\tau-e_{nn}[-1]
\end{matrix}\ts\right|=
\big|\tau+E[-1]\big|_{ii,nn}\cdot\big({-}\tau-e_{nn}[-1]\big),
\een
where $\big|\tau+E[-1]\big|_{ii,nn}$ denotes
the column-determinant of the matrix obtained from
$\tau+E[-1]$ by deleting rows and columns $i$ and $n$.
Since the element $e_{nn}[0]$ annihilates the vacuum vector
of $V(\gl_n)$, the expression for
the determinant in \eqref{firde} simplifies to
\ben
\sum_{k=i+1}^{n-1}\left|\ts\ts\begin{matrix}
\cdots&0&\cdots\\
\cdots&\cdots&\cdots\\
\cdots&\cdots&\cdots\\
\cdots&e_{nk}[-1]&\cdots
\end{matrix}\ts\ts\right|
-\big|\tau+E[-1]\big|_{ii}
+\big|\tau+E[-1]\big|_{ii,nn}\ts\big({-}\tau-e_{nn}[-1]\big).
\een

Applying now Lemma~\ref{lem:comcol} in a similar way
to the determinant in \eqref{ifix}
with $e_{nn}[0]$ replaced with $0$,
we get
\begin{align}
\label{secde}
&\left|\ts\ts\begin{matrix}
\cdots&0&\cdots\\
\cdots&\cdots&\cdots\\
\cdots&0&\cdots\\
\cdots&\cdots&\cdots\\
\cdots&e_{ni}[0\tss]&\cdots
\end{matrix}\ts\ts\right|=(-1)^{n-i}
\left|\ts\ts\begin{matrix}
\cdots&\cdots&0\\
\cdots&\cdots&\cdots\\
\cdots&\cdots&0\\
\cdots&\cdots&\cdots\\
\cdots&\cdots&e_{ni}[0\tss]
\end{matrix}\ts\ts\right|\\
{}+(-1)^{n-i}&\sum_{k=i+1}^{n-1}\left|\ts\ts\begin{matrix}
\cdots&0&\cdots\\
\cdots&\cdots&\cdots\\
\cdots&e_{nk}[-1]&\cdots\\
\cdots&\cdots&\cdots\\
\cdots&0&\cdots
\end{matrix}\ts\ts\right|
+(-1)^{n-i}\left|\ts\ts\begin{matrix}
\cdots&-e_{1\tss i}[-1]\\
\cdots&\cdots\\
\cdots&e_{nn}[-1]-e_{i\tss i}[-1]\\
\cdots&\cdots\\
\cdots&-e_{n-1\tss i}[-1]
\end{matrix}\ts\ts\right|,
\non
\end{align}
where, as before, the dots indicate the same
entries as in the column-determinant
$\cdet\big(\tau+E[-1]\big)$,
the determinants in the second line are of
the size $(n-1)\times(n-1)$ with row $n$ and column $i$ deleted,
and the entry $e_{nk}[-1]$ in the summation term
occur in row $i$ and column $k-1$. Since
the element $e_{ni}[0]$ annihilates the vacuum vector
of $V(\gl_n)$, the first determinant on the right hand side
of \eqref{secde} vanishes. In order to transform the
remaining combination, add and subtract the $(n-1)\times(n-1)$
column-determinant
\ben
(-1)^{n-i}\left|\ts\ts\begin{matrix}
\cdots&\cdots&0\\
\cdots&\cdots&\cdots\\
\cdots&\cdots&-\tau\\
\cdots&\cdots&\cdots\\
\cdots&\cdots&0
\end{matrix}\ts\ts\right|,
\een
obtained from $\cdet(\tau+E[-1])$ by deleting
row $n$ and column $i$ and replacing the last column
as indicated, with $-\tau$ in the row $i$.
After combining this difference with the last term in the expansion
\eqref{secde} we obtain the difference of two
$(n-1)\times(n-1)$
column-determinants multiplied by $(-1)^{n-i}$, where one of
them is obtained from $\cdet(\tau+E[-1])$ by deleting
row $n$ and moving column $i$ to replace the last column.
By Lemma~\ref{lem:tauema}, the corresponding $(n-1)\times(n-1)$
matrix is a Manin matrix. Hence,
taking the signs
into account we conclude that the determinant
in question equals the minor $\big|\tau+E[-1]\big|_{nn}$.
For the remaining determinants we use the general property
of column-determinants allowing us to permute rows.
Moving row $i$ in each of the determinants down to the last position
and taking signs into account, we find that the
determinant in \eqref{secde} equals
\ben
-\sum_{k=i+1}^{n-1}\left|\ts\ts\begin{matrix}
\cdots&0&\cdots\\
\cdots&\cdots&\cdots\\
\cdots&\cdots&\cdots\\
\cdots&e_{nk}[-1]&\cdots
\end{matrix}\ts\ts\right|
+\big|\tau+E[-1]\big|_{nn}
-\big|\tau+E[-1]\big|_{ii,nn}\ts\big({-}\tau-e_{nn}[-1]\big).
\een
Combining this with the
expression for the determinant \eqref{firde} derived above,
we complete the verification of \eqref{ifix}.

Recalling now the expression \eqref{sumdetsi},
we arrive at the relation
\ben
e_{nn}[1]\ts \cdet\big(\tau+E[-1]\big)=\frac{K+n}{n}\ts
\Big((n-1)\ts \big|\tau+E[-1]\big|_{nn}-\sum_{i=1}^{n-1}
\big|\tau+E[-1]\big|_{ii}\Big)
\een
thus completing the proof of \eqref{annih}.

Finally, under the embedding
$\Sr(\gl_n)\hookrightarrow \Sr(\g_-)$ defined by
$e_{ij}\mapsto e_{ij}[-1]$,
the elements $\overline S_1,\dots,\overline S_n$
coincide with the images
of the respective coefficients of the characteristic polynomial
$\det[\de_{ij}\ts u+e_{ij}]$ which
are algebraically independent and
generate the algebra
of invariants $\Sr(\gl_n)^{\gl_n}$.
Therefore, $S_1,\dots,S_n$
form a complete set of Segal--Sugawara vectors
in $V(\gl_n)$.

\section{Center of the local completion and commutative subalgebras}
\label{sec:clc}
\setcounter{equation}{0}

Applying Theorems~\ref{thm:sugaw} and \ref{thm:traces}, we can get
complete sets of Sugawara fields in an explicit form.
Consider the vertex algebra $V(\gl_n)$ and set
$e_{ij}(z)=Y(e_{ij}[-1],z)$ so that
\ben
e_{ij}(z)=\sum_{r\in\ZZ}e_{ij}[r]\ts z^{-r-1},
\qquad i,j=1,\dots,n.
\een
Introduce the
$n\times n$ matrix
$\di_z+E(z)$ by
\ben
\di_z+E(z)=\left[\begin{matrix}
\di_z+e_{11}(z)&e_{12}(z)&\dots&e_{1n}(z)\\
e_{21}(z)&\di_z+e_{22}(z)&\dots&e_{2n}(z)\\
\vdots&\vdots& \ddots&\vdots     \\
e_{n1}(z)&e_{n2}(z)&\dots&\di_z+e_{nn}(z)
                \end{matrix}\right]
\een
and expand its normally ordered column-determinant
\beql{defsf}
:\cdet(\di_z+E(z)):{}=\di_z^{\tss n}+S_1(z)\ts\di_z^{\tss n-1} +\dots
+S_{n-1}(z)\ts\di_z +S_n(z).
\eeq
Equivalently, the fields $S_l(z)$ are given by
$S_l(z)=Y(S_l,z)$, where the elements $S_l$ are defined
in \eqref{defsl}.

\bex\label{ex:cenlc}
For $n=2$ we have
\ben
\bal
:\cdet\left[\begin{matrix}
\di_z+e_{11}(z)&e_{12}(z)\\
e_{21}(z)&\di_z+e_{22}(z)
\end{matrix}\right]:{}
{}&{}={}:\big(\di_z+e_{11}(z)\big)
\big(\di_z+e_{22}(z)\big)
-e_{21}(z)\ts e_{12}(z):
\eal
\een
so that
\ben
S_1(z)=e_{11}(z)+e_{22}(z),\qquad S_2(z)={}
:e_{11}(z)\ts e_{22}(z)-e_{21}(z)\ts e_{12}(z):{}+{}e'_{22}(z).
\een
\vskip-1.2\baselineskip
\qed
\eex

Introduce the fields
$T_{kl}(z)=Y(T_{kl},z)$ corresponding to the
Segal--Sugawara vectors $T_{kl}$.
More explicitly, they
can be defined by the expansion of the normally ordered trace
\beql{deftf}
:\tr\big(\di_z+E(z)\big)^k:{}
= T_{k\tss 0}(z)\ts\di_z^{\tss k}+T_{k1}(z)\ts\di_z^{\tss k-1} +\dots
+T_{kk}(z),\qquad k\geqslant 0.
\eeq

Note that the Sugawara operators associated with
the vectors $T_{kl}$ exhibit some similarity
with the families of operators constructed in
\cite{gw:ho} and \cite{h:so},
although the exact relation with those families is unclear.

Recall that $\z(\wh\gl_n)$ is the algebra
of Segal--Sugawara vectors in
$V(\gl_n)$ which can also be regarded as the
center of the vertex algebra $V_{-n}(\gl_n)$.
The center of the local completion $\U_{-n}(\wh\gl_n)_{\loc}$
is the vector subspace $\Zgot(\wh\gl_n)$ which consists
of the elements commuting with the action of $\wh\gl_n$.
It was proved in \cite{ff:ak} that any element
of $\Zgot(\wh\gl_n)$ is a Fourier component
of a field corresponding to an element of $\z(\wh\gl_n)$.
By Theorems~\ref{thm:sugaw} and \ref{thm:traces}, the fields $Y(S,z)$
with $S\in\z(\wh\gl_n)$ can be interpreted
as differential polynomials either
in the fields $S_1(z),\dots,S_n(z)$
or in the fields $T_{11}(z),\dots,T_{nn}(z)$
with normally ordered products; see \eqref{defsf} and \eqref{deftf}.
Hence, we obtain
two explicit descriptions of $\Zgot(\wh\gl_n)$
formulated below.
The first of them was originally given in \cite{ct:qs},
and the arguments there rely on the results of
\cite{ff:ak}, \cite{r:uh}
and \cite{t:qg}.

\bco\label{cor:cenb} The center
$\Zgot(\wh\gl_n)$ of the local completion
$\U_{-n}(\wh\gl_n)_{\loc}$
consists of the Fourier coefficients of all differential
polynomials in either family of
the fields $S_1(z),\dots,S_n(z)$
or $T_{11}(z),\dots,T_{nn}(z)$.
\qed
\eco

By the vacuum axiom of a vertex algebra,
the application of the fields $S_l(z)$ and $T_{kl}(z)$
to the vacuum vector $1$ of $V_{-n}(\gl_n)$
yields power series in $z$ which we denote respectively by
\beql{posi}
S_l(z)_+=\sum_{r<0}S^+_{l,(r)}z^{-r-1}\Fand
T_{kl}(z)_+=\sum_{r<0}T^+_{kl,(r)}z^{-r-1}.
\eeq
Since
the Fourier coefficients of the fields $S_l(z)$ and $T_{kl}(z)$
belong to $\Zgot(\wh\gl_n)$, all coefficients of the series
\eqref{posi} belong to the center $\z(\wh\gl_n)$
of the vertex algebra $V_{-n}(\gl_n)$.
By \eqref{defsf} and \eqref{deftf}, the series
\eqref{posi} can be written in an explicit form
with the use of the matrix
$E(z)_+=[e_{ij}(z)_+]$. We have
\ben
\cdet(\di_z+E(z)_+)=\di_z^{\tss n}+ S_1(z)_+\ts\di_z^{\tss n-1}+\dots
+S_{n-1}(z)_+\ts\di_z +S_n(z)_+
\een
and
\ben
\tr\big(\di_z+E(z)_+\big)^k
= T_{k\tss 0}(z)_+\ts\di_z^{\tss k}
+ T_{k1}(z)_+\ts\di_z^{\tss k-1}+\dots
+T_{kk}(z)_+.
\een
We arrive at the following result, whose first part
dealing with the commuting family of the elements $S^+_{l,(r)}$
goes back to the
original work \cite{t:qg};
see also \cite{cf:mm}, \cite{ct:qs}.

\bco\label{cor:coomsa}
The elements of each of the families
\ben
S^+_{l,(r)}\quad\text{with}\quad l=1,\dots,n\fand r< 0,
\qquad
T^+_{kl,(r)}\quad\text{with}\quad
0\leqslant l\leqslant k\fand r< 0,
\een
belong to $\z(\wh\gl_n)$. Moreover,
$\z(\wh\gl_n)$ is the algebra of polynomials
\ben
\z(\wh\gl_n)=\CC[S^+_{l,(r)}\ |\ l=1,\dots,n,\ \ r< 0]
=
\CC[T^+_{l\tss l,(r)}\ |\
l=1,\dots,n;\ r< 0\tss].
\een
\eco

\bpf
We only need to prove the
algebraic independence of the families of
the elements $S^+_{l,(r)}$ and $T^+_{l\tss l,(r)}$.
This follows by comparing their highest degree
components in the graded algebra
$\text{gr}\ts\U(\g_-)\cong \Sr(\g_-)$
with those of the elements $T^r\tss S_l$ and $T^r\tss T_{l\tss l}$,
respectively; see Theorems~\ref{thm:sugaw} and \ref{thm:traces}.
In the general case
such relationship between two families of generators
of $\Sr(\g_-)$ was pointed out in \cite{r:uh}.
\epf

Note that in our approach,
the pairwise commutativity of the elements $S^+_{l,(r)}$ in
$\U(\g_-)$
is a consequence of
the fact that the $S^+_{l,(r)}$ are Segal--Sugawara vectors.
The higher Gaudin Hamiltonians can be obtained by taking
the images of the elements of the commutative subalgebra
$\z(\wh\gl_n)$ of $\U(\g_-)$
in the algebras $\U(\gl_n)^{\ot m}$ under certain
evaluation homomorphisms; see e.g. \cite{ffr:gm}, \cite{t:qg}
for details.

\section{Eigenvalues in the Wakimoto modules}
\label{sec:ewm}
\setcounter{equation}{0}

As another application of the explicit
formulas of Theorems~\ref{thm:sugaw} and \ref{thm:traces}, we recover
a description of the classical $\Wc$-algebra
$\Wc(\gl_n)$. This allows one to calculate the eigenvalues
of the elements of the center $\Zgot(\wh\gl_n)$
in the Wakimoto modules of critical level. We start by recalling
the main steps in the construction of these modules
following \cite{f:wm}.

Consider the {\it Weyl algebra\/}
$\Ac(\gl_n)$ generated by the elements
$a_{ij}[r]$ with $r\in\ZZ$, $i,j=1,\dots,n$ and $i\ne j$
and the defining relations
\ben
\big[a_{ij}[r],a_{kl}[s]\tss\big]
=\de_{kj}\tss\de_{i\tss l}\tss\de_{r,-s}
\qquad\text{for}\quad i<j,
\een
whereas all other pairs of the generators commute.
The Fock representation $M(\gl_n)$ of $\Ac(\gl_n)$
is generated by a vector $|0\rangle$
such that for $i<j$ we have
\ben
a_{ij}[r]|0\rangle=0,\quad r\geqslant 0\Fand
a_{ji}[r]|0\rangle=0,\quad r>0.
\een
The vector space $M(\gl_n)$ carries a vertex algebra
structure. In particular,
$|0\rangle$ is the vacuum vector, and
for $i<j$ we have
\ben
Y(a_{ij}[-1]\ts|0\rangle,z)=\sum_{r\in\ZZ} a_{ij}[r]\tss z^{-r-1}
\Fand
Y(a_{ji}[0]\ts|0\rangle,z)=\sum_{r\in\ZZ} a_{ji}[r]\tss z^{-r}.
\een
We will denote these series by $a_{ij}(z)$ and $a_{ji}(z)$,
respectively.

Take an $n$-tuple
\ben
\chi(t)=\big(\chi_1(t),\dots,\chi_n(t)\big),\qquad
\chi_i(t)=\sum_{r\in\ZZ}\chi_i[r]\tss t^{-r-1}\in\CC((t)),
\een
where $\CC((t))$ denotes the algebra of formal Laurent series in $t$
containing only a finite number of negative powers of $t$.
The following formulas define a
$\wh\gl_n$-module structure on the vector space $M(\gl_n)$:
\beql{wak}
\bal
e_{i,i+1}(z)&\mapsto a_{i,i+1}(z)+\sum_{k<l}
:P^{\tss i}_{kl}\ts a_{kl}(z):\\
e_{ii}(z)&\mapsto \sum_{k<l}
d_{kl}^{\ts i}\ts:a_{l\tss k}(z)\tss a_{kl}(z):{}+\chi_i(z)\\
e_{i+1,i}(z)&\mapsto \sum_{k<l}
:Q^{\tss i}_{kl}\ts a_{kl}(z):+\ts c_i\ts\di_z\ts a_{i+1,i}(z)
+\big(\chi_i(z)-\chi_{i+1}(z)\big)\ts a_{i+1,i}(z),
\eal
\eeq
where $P^{\tss i}_{kl}$ and $Q^{\tss i}_{kl}$ are certain polynomial
expressions in the $a_{pq}(z)$ with $p>q$, while $d_{kl}^{\ts i}$
and $c_i$ are certain coefficients. The precise formulas
are obtained by applying the formulas
\cite[(4.6)--(4.8)]{f:wm} to the Lie algebra
$\wh\sll_n$ and extending the action
to $\wh\gl_n$ by setting
\ben
e_{11}(z)+\dots+e_{nn}(z)\mapsto \chi_1(z)+\dots+\chi_n(z).
\een
These are the {\it Wakimoto modules of critical level\/},
denoted by $W_{\chi(t)}$.

\bex\label{ex:waktwo}
For $n=2$ the explicit formulas are
\beql{waktwo}
\bal
e_{12}(z)&\mapsto a_{12}(z)\\
e_{11}(z)&\mapsto
\ts{}-{}:a_{21}(z)\tss a_{12}(z):{}+\chi_1(z)\\
e_{22}(z)&\mapsto
{}:a_{21}(z)\tss a_{12}(z):{}+\chi_2(z)\\
e_{21}(z)&\mapsto
\ts{}-{}:a_{21}(z)^2\tss a_{12}(z):{}
-2\ts\di_z\ts a_{21}(z)
+\big(\chi_1(z)-\chi_2(z)\big)\ts a_{21}(z).
\eal
\eeq
\eex

The key fact leading to the construction of the Wakimoto modules
is the existence of the vertex algebra homomorphism
\beql{rhohom}
\rho:V_{-n}(\gl_n)\to M(\gl_n)\ot\pi_0,
\eeq
where $\pi_0$ is
the algebra of polynomials
\ben
\pi_0=\CC[b_i[r]\ |\ i=1,\dots,n;\ r<0\tss]
\een
in the variables $b_i[r]$ which we consider as
a commutative vertex algebra.
The translation operator on $\pi_0$ is defined by
\ben
T\ts 1=0,\qquad \big[T,b_i[r]\tss\big]=-r\ts b_i[r-1].
\een
The homomorphism \eqref{rhohom} is defined by the
formulas \eqref{wak}, where the series
$\chi_i(z)$ should be respectively replaced
by
\ben
b_i(z)=\sum_{r<0}b_i[r]\ts z^{-r-1}.
\een
The image of the center $\z(\wh\gl_n)$ of the vertex algebra
$V_{-n}(\gl_n)$ under the homomorphism $\rho$
is contained
in $\pi_0\cong 1\ot\pi_0$.
Moreover, this image coincides with
the {\it classical $\Wc$-algebra\/}
$\Wc(\gl_n)$; see \cite[Sec.~9]{f:wm}.
In order to recall the definition of
$\Wc(\gl_n)$, introduce the operators
\ben
Q_i:\pi_0\to\pi_0,\qquad i=1,\dots,n-1,
\een
by the formulas
\ben
Q_i=\sum_{r=0}^{\infty}\ts
\sum_{\la\tss\vdash\tss r}\frac{{\mathbf b}_i(\la)}{z_{\la}}\ts
\Big(\frac{\di}{\di\tss b_i[-r-1]}
-\frac{\di}{\di\tss b_{i+1}[-r-1]}\Big),
\een
where the following notation was used.
The second sum is taken over all partitions
$\la=(\la_1,\dots,\la_p)$ of $r$, so that
$\la_1\geqslant\dots\geqslant\la_p>0$ and
$\la_1+\dots+\la_p=r$.
Furthermore,
\ben
{\mathbf b}_i(\la)=\big(b_i[-\la_1]-b_{i+1}[-\la_1]\big)\dots
\big(b_i[-\la_p]-b_{i+1}[-\la_p]\big)
\een
and
\ben
z_{\la}=1^{m_1}m_1!\ts 2^{m_2}m_2!\dots r^{m_r}m_r!,
\een
where $m_k$ denotes the multiplicity
of $k$ in $\la$ for each $k$. By definition, $\Wc(\gl_n)$
consists of the elements of $\pi_0$, annihilated
by all operators $Q_i$,
\ben
\Wc(\gl_n)=\bigcap_{\phantom{-1}1\leqslant\tss i\tss\leqslant n-1}
\text{Ker}\ts Q_i.
\een
Applying Theorems~\ref{thm:sugaw} and \ref{thm:traces},
we recover an explicit description of the algebra $\Wc(\gl_n)$; see
e.g. \cite{d:lc}, \cite{ff:im}. In order to formulate the result,
consider the extended algebra, which is isomorphic to
$\pi_0\ot\CC[\tau]$, as a vector space,
with the relations
\ben
\big[\tau,b_i[r]\tss\big]=-r\ts b_i[r-1] .
\een

\bco\label{cor:walg}
The image of the column-determinant $\cdet(\tau+E[-1])$
under the homomorphism \eqref{rhohom}
is given by
\ben
\rho:\cdet(\tau+E[-1])\mapsto \big(\tau+b_n[-1]\big)\cdots
\big(\tau+b_1[-1]\big).
\een
Hence, if the elements $B_i$ are defined by
\ben
\big(\tau+b_n[-1]\big)\cdots\big(\tau+b_1[-1]\big)
=\tau^n+ B_1\ts\tau^{n-1}+\dots+B_n,
\een
then $\Wc(\gl_n)$ is the algebra of polynomials
in the variables $T^{\tss r}B_i$
with $i=1,\dots,n$ and $r\geqslant 0$.
\eco

\bpf
Using Lemma~\ref{lem:tauema} and the properties
of the Manin matrices, write the rows and columns
of the matrix $\tau+E[-1]$ in their reverse orders
so that the column determinant will take the form
\beql{detexp}
\cdet(\tau+E[-1])=\sum_{\si}\sgn\si\cdot
\big(\de_{\si(n)n}\tss\tau+e_{\si(n)n}[-1]\big)\dots
\big(\de_{\si(1)1}\tss\tau+e_{\si(1)1}[-1]\big).
\eeq
Now we apply the explicit description of the homomorphism
\eqref{rhohom} given in \cite{f:wm}.
Given that
the image of $\z(\wh\gl_n)$
under this homomorphism is contained
in $\pi_0$, we only need to keep the terms
containing the elements $b_i[r]$. The images of
the generators $e_{ij}[r]$ with $i<j$ and $r<0$ do not
depend on the $b_i[r]$ and
the images of the summands in \eqref{detexp} with $\si(n)<n$
do not contribute to the image of
the column determinant. Hence, the nonzero contribution
is only provided by the summands in \eqref{detexp} with $\si(n)=n$.
Considering now $\si(n-1)$ and arguing by induction, we
conclude that the only summand in \eqref{detexp}
providing a nonzero contribution to the image
corresponds to the identity permutation $\si$.
This proves the first part.
The second statement follows from \cite[Theorem~9.6]{f:wm}.
\epf

Note that an alternative calculation of the image of
$\cdet(\tau+E[-1])$ under
the homomorphism \eqref{rhohom}
can be performed with the use
of the original expression provided by \eqref{cdetdef}
and the application of
\cite[Theorem~11.3]{f:wm}.

\bex\label{ex:genwa}
If $n=2$ then
\ben
\bal
B_1&=b_{1}[-1]+b_{2}[-1],\\
B_2&=b_{1}[-1]\ts b_{2}[-1]+b_{1}[-2].
\eal
\een
If $n=3$ then
\ben
\bal
B_1&=b_{1}[-1]+b_{2}[-1]+b_{3}[-1],\\
B_2&=b_{1}[-1]\ts b_{2}[-1]+b_{1}[-1]\ts b_{3}[-1]
+b_{2}[-1]\ts b_{3}[-1]+2\ts b_{1}[-2]+b_{2}[-2],\\
B_3&=b_{1}[-1]\ts b_{2}[-1]\ts b_{3}[-1]+b_{1}[-2]\ts b_{2}[-1]\\
{}&\qquad\qquad\qquad\qquad
+b_{1}[-2]\ts b_{3}[-1]+b_{1}[-1]\ts b_{2}[-2]+2\ts b_{1}[-3].
\eal
\een
\vskip-0.8\baselineskip
\qed
\eex

Applying the identity \eqref{newton}
and Corollary~\ref{cor:walg}, we can also
calculate the images of the Segal--Sugawara vectors $T_{kl}$.

\bco\label{cor:traim}
The images of the elements $\tr(\tau+E[-1])^k$
under the homomorphism \eqref{rhohom}
are found from the formula
\begin{multline}
\rho:\sum_{k=0}^{\infty}t^k\ts \tr(\tau+E[-1])^k
\mapsto
\sum_{i=1}^n\Big(1-t\tss\big(\tau+b_1[-1]\big)\Big)^{-1}\cdots
\Big(1-t\tss\big(\tau+b_{i}[-1]\big)\Big)^{-1}\\
{}\times{}\Big(1-t\tss\big(\tau+b_{i-1}[-1]\big)\Big)\cdots
\Big(1-t\tss\big(\tau+b_1[-1]\big)\Big),
\non
\end{multline}
where $t$ is a complex variable.
\eco

\bpf
The formula follows from \eqref{newton}
by calculating the derivative on the left hand side
and then replacing $u$ by $-t^{-1}$.
\epf

The elements of the center $\Zgot(\wh\gl_n)$ of
$\U_{-n}(\wh\gl_n)_{\loc}$ act on the
Wakimoto modules $W_{\chi(t)}$
as multiplications by scalars which can be calculated
by using Corollaries~\ref{cor:walg} and \ref{cor:traim}.
Recall the description of $\Zgot(\wh\gl_n)$
provided by Corollary~\ref{cor:cenb}.

\bco\label{cor:eigw}
The coefficients of $:\cdet(\di_z+E(z)):$
and $:\tr\big(\di_z+E(z)\big)^k:$
act on $W_{\chi(t)}$
as multiplications by scalars found from the
respective formulas
\ben
:\cdet(\di_z+E(z)):{}\mapsto
\big(\di_z+\chi_n(z)\big)\dots \big(\di_z+\chi_1(z)\big)
\een
and
\begin{multline}
\sum_{k=0}^{\infty}t^k\ts
:\tr\big(\di_z+E(z)\big)^k:{}
\mapsto
\sum_{i=1}^n\Big(1-t\tss\big(\di_z+\chi_1(z)\big)\Big)^{-1}\cdots
\Big(1-t\tss\big(\di_z+\chi_i(z)\big)\Big)^{-1}\\
{}\times{}\Big(1-t\tss\big(\di_z+\chi_{i-1}(z)\big)\Big)\cdots
\Big(1-t\tss\big(\di_z+\chi_1(z)\big)\Big).
\non
\end{multline}
\eco

\section*{Acknowledgments}

We acknowledge the support of the Australian Research
Council. The first author is grateful to the University of
Sydney for hospitality during his visit.
His research is also supported
by the grant of Support for the
Scientific Schools 8004.2006.2,  RFBR grant 08-02-00287a,
the ANR grant GIMP (Geometry and Integrability in
Mathematics and Physics) and the grant MK-5056.2007.1.

\end{document}